\documentclass[11pt]{amsart}
\usepackage{graphicx}
\usepackage{setspace}
\usepackage{eufrak}
\usepackage{float}
\usepackage{etex}
\usepackage[all]{xy}
\usepackage{amssymb, amsmath}
\usepackage{ascmac}

\newcommand{\Z}{\mathbb{Z}}
\newcommand{\N}{\mathbb{N}}

\newcommand{\R}{\mathbb{R}}

\newcommand{\Int}{\mbox{\rm Int}\, }

\def\t{\noindent}
\def\rd{\partial}

\def\S{\mathcal{S}}

\def\bu{\mbox{\boldmath $u$}}

\newtheorem{thm}{\bf Theorem}[section]
\newtheorem{cor}[thm]{\bf Corollary} 
\newtheorem{lem}[thm]{\bf Lemma} 
\newtheorem{prop}[thm]{\bf Proposition} 
\newtheorem{definition}[thm]{\bf Definition} 
 
\newtheorem{rem}[thm]{\bf Remark}

\begin{document}
\setlength{\baselineskip}{16pt}

\title[$C^1$-triangulations of semialgebraic sets]
{$C^1$-triangulations of semialgebraic sets} 
\author[T.~Ohmoto]{Toru Ohmoto}
\address[T.~Ohmoto]{Department of Mathematics, 
Faculty of Science,  Hokkaido University,
Sapporo 060-0810, Japan}
\email{ohmoto@math.sci.hokudai.ac.jp}
\author[M.~Shiota]{Masahiro Shiota}
\address[M.~Shiota]{Graduate School of Mathematics, 
Nagoya University,Nagoya, 464-8602, Japan}
\email{shiota@math.nagoya-u.ac.jp}
%
%
%
\maketitle
\begin{abstract}
We show that every semialgebraic set admits a semialgebraic triangulation such
that each closed simplex is  $C^1$ differentiable. 
As an application, we give a straightforward definition of 
 the integration $\int_X \omega$  over a compact semialgebraic subset $X$ 
 of a differential form $\omega$ on an ambient semialgebraic manifold. 
 This provides a significant simplification of 
 the theory of semialgebraic singular chains and integrations 
 without using geometric measure theory. 
 Our results hold over every (possibly non-archimedian) real closed field. 
\end{abstract}
%
%

\section{Introduction}

In this paper, we solve a fundamental question 
on the $C^1$-regularity of triangulations for semialgebraic sets. 
Also an elementary application to the integration of differential forms is discussed, 
that would be useful, e.g. for the de Rham homotopy theory of semialgebraic sets introduced by  
Kontsevich-Soibelman \cite{KS}. 
For simplicity, we mainly work on the category of semialgebraic sets and maps 
over the real number field $\R$, 
although our argument goes well in a more general category 
over a general real closed field, e.g. 
an {\it o-minimal} category, the category of {\it subanalytic sets and maps}   
or more generally an $\frak X$-category with Axiom (v) in \cite[II, p.97]{S}. 

To begin with, we explain a convention of terminologies 
which will be used throughout the present paper. 
Let $X\subset\R^m$ and $Y\subset\R^n$ be semialgebraic sets,  
{\it i.e.}, subsets defined by finitely many polynomial equations, inequalities and Boolean operations; 
a semialgebraic map $X\to Y$ is a map whose graph is a semialgebraic set in $\R^m \times \R^n$. 
A subset $X \subset \R^m$ is 
a {\it locally semialgebraic set}  if for any compact $m$-disk $D$ in $\R^m$,  
$X \cap D$ is semialgebraic; 
a map $X\to Y$ is called {\it locally semialgebraic} 
if for any compact semialgebraic subset $X'$ of $X$, the restriction to $X'$ is a semialgebraic map. 
Note that locally semialgebraic does not imply semialgebraic. 
Because we systematically deal with possibly non-compact spaces, 
{\it we use the notion semialgebraic to mean locally semialgebraic in the above sense}, 
unless specifically mentioned. 
Let $k$ be an integer. 
A {\it semialgebraic $C^k$ map} $X\to Y$ is the restriction of a semialgebraic map of class $C^k$ from some semialgebraic open neighborhood of $X$ in $\R^m$ to $\R^n$; 
in fact, we can extend it to $\R^m\to \R^n$ by multiplying 
a semialgebraic $C^k$ function $\phi$ on $\R^m$ such that $\phi=1$ on $X$ 
and $\phi=0$ outside of a smaller closed neighborhood (cf.  \cite{S}). 

We write (half-)open intervals like as 
$(0,1)=\{\, t\in \R\,  | \, 0<t<1 \, \}$ 
and $[0,1)=\{\, t\in \R\,  | \, 0\le t<1 \, \}$ etc. 

\subsection{$C^1$ triangulation} 
In the present paper,  
a {\it semialgebraic triangulation} of a locally closed semialgebraic set $X$ means    
the pair of a locally finite simplicial complex $K$ 
and a semialgebraic homeomorphism $f:|K|\to X$, 
where $|K|$ is the polyhedron in some Euclidean space $\R^\ell$ 
obtained as the geometric realization of $K$. 
It is well known that 
such a semialgebraic triangulation always exists (see Theorem \ref{triangulation} below); 
moreover, 
it is possible to take $(K, f)$ so that 
$f|_{{\tiny \Int}(\sigma)}$ is a $C^k$ diffeomorphism for each $\sigma \in K$, 
where $\Int (\sigma)$ denotes the interior of $\sigma$ and $k\ge 1$. 
However we here do not require any regularity of  $C^0$ manifolds $f(\Int(\sigma))$ in 
our definition of semialgebraic triangulation. 
Instead, 
we are concerned with 
the differentiability (the regularity) of $f(\sigma)$ along the boundary $\rd \sigma$.  
Our main theorem is the following: 

\begin{thm}\label{thm1} For a locally closed semialgebraic set $X$, 
there exists a semialgebraic triangulation $(K, f)$ 
so that the map $f:|K|\to X$ is of class $C^1$. 
\end{thm}

We call $(K, f)$ in the theorem 
a {\it semialgebraic triangulation of $X$ with a $C^1$ realization}. 
Then, 
it follows that for each $\sigma \in K$, 
$f(\sigma)$ is semialgebraically homeomorphic to the closed simplex $\sigma \subset |K|$ 
 and the differential $d(f|_\sigma)$ exists and is continuous on $\sigma$ 
 (precisely saying, 
 $f$ is of class $C^1$ on a neighborhood $U_\sigma$ of $\sigma$ in $\R^\ell$, and hence, 
 taking the affine subspace $H$ of the same dimension containing $\sigma$, 
 $f|_{H\cap U_\sigma}$ is also $C^1$). 
 However,  as a cost for that, 
$d(f|_\sigma)$ may drop the rank along the boundary $\rd \sigma$ and 
some smaller dimensional semialgebraic subset of $\Int(\sigma)$. 

The proof is elementary, 
only using standard techniques such as triangulation, tube and the curve selection lemma 
(some basic facts are summarized in Section 2), 
so it would be readable even for non-experts. 
The proof is done in the case of semialgebraic sets over $\R$, 
but it is presented in such a way that it is valid in a more general context 
of $\frak X$-category over a real closed field $R$ in the sense of \cite{S}; 
in particular, it works in the category of subanalytic sets and maps.  
Furthermore, our argument fits with 
{\it piecewise algebraic spaces} introduced in \cite{KS} 
(Remark \ref{PA}). 

There remain some open questions. For instance, 
we do not know 
whether or not $f$ can be of class $C^k$ ($2\le k\le \infty$ or $k=\omega$). 
As known, any Whitney stratified set is triangulable, 
so it would be natural to ask 
if the realization map $f$ can be $C^1$,  see  e.g. \cite{MT}. 

\subsection{Differential forms and integration} 
Let $M$ be a semialgebraic manifold 
and $X$ a compact semialgebraic subset of $M$ ($p=\dim X$). 
Take a semialgebraic triangulation $(K, f)$ of $X$ with a $C^1$ realization. 
Here $K$ is a finite simplicial complex. 
For a differential $p$-form $\omega \in \Omega^p(M)$ and a $p$-simplex $\sigma \in K$, 
 the integral of $\omega$ on $\sigma$ is defined in an obvious way by 
$$\int_{f(\sigma)} \omega := \int_\sigma (f|_\sigma)^* \omega$$
because the map $f|_\sigma: \sigma \to M$ is of class $C^1$  
(the integral is often written by $\int_\sigma \omega$ for short). 
Suppose that $X$ has the $\Z$-fundamental class. 
Put 
$$\int_X \omega := \sum \int_{f(\sigma)} \omega$$ 
where the sum runs over all $p$-simplices $\sigma$ in $K$ with suitable orientation 
so that they form the fundamental cycle. 
It is almost obvious that  the value $\int_X \sigma$ 
does not depend on the choice of $(K, f)$ of $X$ (Theorem \ref{thm3}). 
Consequently, 
the theory of integration of differential forms over semialgebraic sets 
can be treated in entirely the same manner as in the case of manifolds  -- 
that simplifies the treatment of semialgebraic chains and integrations, 
without using geometric measure theory due to Hardt et al 
\cite{Federer, Hardt, HLTV, KS} (see Remarks \ref{current}, \ref{CDGA}). 
In fact, a traditional way using geometric currents involves the limiting process: 
$\int_\sigma \omega$ is defined by 
the integral of $\omega$ over the interior $\Int (\sigma)$ which absolutely converges 
(that is shown by the fact that the induced volume is finite). 
On the other hand, our definition uses only  
the {\it definite integral} of a continuous function on  a compact domain $\sigma$. 

We remark that our definition of 
$\int_X \omega$ works over 
not only the real number field $\R$ 
but also any possibly non-archimedian real closed field $R$. 
Then the value of the integral does not take values in the same field $R$ but 
in a larger saturated field (see \cite{MS}). 
In particular, the volume of a bounded definable set in $R^m$ is bounded within $R$ 
(Remark \ref{o_minimal}).

\subsection*{Acknowledgement} 
This work was originally motivated by a question of Tatsuo Suwa 
about the definition of $\int_{X} \omega$ using triangulation, 
which was needed for rigorously establishing the Alexander duality 
for an embedded analytic subvariety 
in the \v{C}ech-de Rham context \cite{Suwa}. 
The authors would like to thank Tatsuo Suwa for guiding them to this problem 
and for useful discussions. 
During preparing the final version of the present paper,  
we were informed about two recent preprints \cite{CP} and \cite{HKT}; 
we also thank G. Comte and M. Hanamura for letting us know of those papers, respectively. 
The first author is  supported by JSPS KAKENHI Grants no.24340007 
and no.15K13452, 
and the second author is supported by JSPS KAKENHI Grant no.26400084.

\section{Preliminary} 
We summarize some basic facts in semialgebraic geometry,  
which are also valid without any change for the {\it locally semialgebraic} setting in our sense. 
As noted in Introduction, we simply say semialgebraic to mean locally semialgebraic, 
unless specifically mentioned. 

\subsection{Triangulation} 
Let $k$ be an integer and $X \subset \R^m$ a semialgebraic set. 
A {\it semialgebraic $C^k$ stratification of $X$} is 
a family $\S=\{X_i\}$ of semialgebraic subsets of $X$ (called {\it strata} of $\S$)  
so that each $X_i$ is a locally closed $C^k$ submanifold of $\R^m$, 
$X=\bigsqcup X_i$ (disjoint), $X_i \cap X_j=\emptyset \; (i\not=j)$,  $\S$ is locally finite, 
and $\S$ satisfies the {\it frontier condition}, 
{\it i.e.}, for each $X_i$, 
the closure $\overline{X_i}$ is a union of $X_i$ and some $X_j \in \S$ of smaller dimension. 
We say that $\mathcal{S}$ is {\it compatible with} a family $\mathcal{A}$ of subsets of $X$ 
 if for each stratum $V \in \mathcal{S}$ and each $A \in \mathcal{A}$, it holds either 
 $V\subset A$ or $V\cap A=\emptyset$. 
A {\it semialgebraic $C^k$ stratification of a semialgebraic $C^0$ map $f:X\to Y$} 
is the pair of semialgebraic $C^k$ stratifications $\mathcal{S}$ of $X$ 
and $\mathcal{T}$ of $Y$ such that for each $V \in \mathcal{S}$, 
$f(V)$ is a stratum of $\mathcal{T}$ 
and $f|_{V}:V\to f(V)$ is a semialgebraic $C^k$ submersion.
We denote it by $f:(X; \mathcal{S})\to (Y; \mathcal{T})$. 
Such a semialgebraic stratification always exists 
(the statement can be more enhanced so that 
$\mathcal{S}$ and $\mathcal{T}$ satisfy the Whitney condition, 
but we do not need it later):

\begin{lem} \label{C1_stratif} 
{\rm ({\bf Canonical stratification of maps} \cite{G}, \cite[I. Lemmas 2.2, 2.6; 2.10, 2.14]{S})} 
Let $X$ be a semialgebraic set and $k\ge 1$. 
Let  $X'$ be a closed semialgebraic subset of $X$, $\mathcal{S}$ 
a semialgebraic $C^k$ stratification of $X-X'$, and $\mathcal{A}$ 
a finite family of semialgebraic subsets of $X'$.
Assume $\dim X'<\dim V$ for any $V \in \mathcal{S}$. 
Then there exists a semialgebraic $C^k$  stratification 
$\mathcal{S}'$ of $X'$ compatible with $\mathcal{A}$ 
such that $\mathcal{S}\sqcup\mathcal{S}'$ is 
a semialgebraic  $C^k$ stratification of $X$.\par
Moreover, let $f:X\to Y$ be a proper semialgebraic $C^0$ map 
such that $f(X')\cap f(X-X')=\emptyset$.
Assume $f:(X-X'; \mathcal{S}) \to (f(X-X'); \mathcal{T})$, 
where $\mathcal{T}$ consists of 
all $f(V)$, $V \in \mathcal{S}$, 
is a semialgebraic $C^k$ stratification of $f|_{X-X'}$.
Then we can choose $\mathcal{S}'$ so that 
$f:\mathcal{S}\sqcup\mathcal{S}'\to\mathcal{T}\sqcup\mathcal{T}'$, 
where $\mathcal{T}'$ consists of 
all $f(V')$, $V' \in \mathcal{S}'$, 
 is a semialgebraic $C^k$ stratification of $f:X\to f(X)$.\end{lem}

\begin{thm}\label{triangulation}  
{\rm ({\bf {\L}ojasiewicz's triangulation theorem} \cite{L}, \cite[II.2.1]{S})}
Let $X$ be a locally closed semialgebraic set in $\R^m$ and 
$\mathcal{A}$ a finite family of semialgebraic subsets of $X$. 
Then there is a pair of $(K, h)$ of a locally finite simplicial complex $K$ 
and a semialgebraic homeomorphism $h: |K| \to X$ 
so that  $\{h(\Int(\sigma)), \; \sigma \in K\}$ is 
a semialgebraic $C^k$ stratification of $X$ 
compatible with $\mathcal{A}$. 
\end{thm}

Those theorems have been extended to more general categories 
(o-minimal, subanalytic, $\frak X$), e.g. see \cite{S}.

\subsection{Curve selection lemma}
We will frequently use the following version of 
the {\it curve selection lemma} ({\it wing lemma}) afterward. 
There is a number of literature; the origin goes back to Whitney and Milnor, and 
a recent one is, e.g. a version in o-minimal category \cite[Lem.1.7]{Loi}. 
For proving its $\frak X$-version, 
it suffices to use the $\frak X$-isotopy lemma \cite[II.6]{S}.  
 
\begin{lem}\label{csl}
{\rm ({\bf Wing lemma})} 
Let $V \subset \R^m$ be a non-empty semialgebraic open subset,  
and $Z$ a semialgebraic subset of $\R^m\times \R^n$ so that 
its closure $\overline{Z}$ contains $V \times \{0\}$ and 
$Z\cap (V \times \{0\}) =\emptyset$. 
Then 
there is a non-empty semialgebraic open subset $U \subset V$ 
which admits 
a semialgebraic $C^0$ map $c: U\times [0, a) \to \R^n$ 
with some $a>0$ so that 
$c(x,0)=0$ for $x \in U$, $(x, c(x,t)) \in Z$ and $|c(x,t)|=t$ for any $t\not=0$,  
and $c$ is of class $C^k$ on $U\times (0, a)$. 
\end{lem}

\section{Proof -- Panel beating} 

\subsection{Panel beating}
We will prove a bit more general statement. 

\begin{thm}\label{thm2}
Let $X$ be a locally closed semialgebraic set with a semialgebraic triangulation $(K, f)$, 
and  $\varphi: X \to Y$ a semialgebraic $C^0$ map. 
Then there exists a semialgebraic homeomorphism $\chi$ of $|K|$ 
such that $\chi$ preserves any $\sigma \in K$ and 
the composed map 
$$
\varphi \circ f \circ \chi: 
|K| \stackrel{\simeq}{\longrightarrow}  |K| \stackrel{\simeq}{\longrightarrow} X \longrightarrow Y 
$$
is of class  $C^1$. 
\end{thm}

Theorem \ref{thm1} is the case of $X=Y$ and $\varphi=id_X$. 

\begin{definition}
{\rm 
We call $\chi$ a {\it panel beating} of the triangulation $(K, f)$ of $X$ 
with respect to $\varphi: X \to Y$.  
}
\end{definition}

Notice that 
in Theorem \ref{thm2}, 
$\chi$ can be chosen so that $\chi$ itself, $f\circ \chi$ and 
$\varphi \circ f\circ \chi$ are of class $C^1$ simultaneously. In fact, we have 

\begin{cor}\label{thm2_cor}
Let $X$ and $(K, f)$ as above. 
Suppose that we are given  semialgebraic  $C^0$ maps 
$\varphi^{(j)}: X \to Y_j\; (1\le j \le s)$. 
Then, there is a common panel beating $\chi$ 
such that each $\varphi^{(j)}\circ f \circ \chi$ is of class $C^1$. 
\end{cor}

\proof 
Apply Theorem \ref{thm2} to the map 
$\varphi=\prod \varphi^{(j)}: X \to Y=\prod Y_j$. \qed

\

Here is an obvious example. 
Consider $\varphi: \R \to \R$ given by $\varphi(u)=\sqrt{|u|}$, 
and decompose $\R$ into $(-\infty, 0)$, $\{0\}$ and $(0, \infty)$. 
Put 
$$\chi_0(u)=\left\{
\begin{array}{rc}
u^4 & (u\ge 0)\\
-u^4 & (u\le 0)
\end{array}\right. 
$$
then it is a panel beating; indeed, $\varphi\circ \chi_0(u)=u^2$. 
At a key step  (Lemma \ref{lemma1}) in the proof of Theorem \ref{thm2}, 
we glue together this kind of mappings $\chi_0$ 
via a semialgebraic version of partition of unity 
to construct the desired semialgebraic homeomorphism $\chi:|K| \to |K|$.

\subsection{Proof of Theorem \ref{thm2}} 
Let $X$ be a locally closed semialgebraic set 
which admits a semialgebraic triangulation $(K_1, f_1)$ 
with $|K_1| \subset \R^m$, and 
$Y$ a semialgebraic set in $\R^p$. 
Let $\varphi: X \to Y$ be a semialgebraic $C^0$ map. 
By the semialgebraic Tietze theorem, e.g. \cite[Prop.2.6.9]{BCR} 
(or using a semialgebraic partition of unity and tubes \cite{S}), 
we see that the map $\varphi \circ f_1$ 
is extended to a semialgebraic $C^0$ map $\R^m \to \R^p$. 
The graph map is denoted by 
$f:\R^m \to \R^n:=\R^m\times \R^p$, 
which is semialgebraic $C^0$ and proper. 
Also we extend $K_1$ to a simplicial decomposition 
(locally closed infinite simplicial complex) $K$ of $\R^m$. 
To prove Theorem \ref{thm2} 
we find a semialgebraic homeomorphism 
$\chi: \R^m \to \R^m$ so that 
$\chi$ preserves any simplices of $K$ and $f\circ \chi$ is of class $C^1$. 

By Lemma \ref{C1_stratif}, 
there is a semialgebraic $C^2$ stratification $\mathcal{S}\to\mathcal{T}$ of $f$,  
such that each simplex in $K$ is a finite union of some strata of $\mathcal{S}$. 
The restriction of $f$ on each stratum is of class $C^2$ by definition, 
so in particular $f$ is $C^1$ on the union of open strata. 
Now we use the downward induction on 
dimension $d\; (\le m-1)$ of the `bad subset' where $f$ is not $C^1$: 
We will construct a panel beating $\chi$ around the subset 
so that $f\circ \chi$ is $C^1$ off a smaller subset of less dimension. 
Clearly, it suffices to construct $\chi$ around each connected stratum of top dimension 
of the bad subset. 
Furthermore, we may assume that 
the stratum is semialgebraically diffeomorphic to $\R^d$ 
 (by subdividing $\mathcal{S}$ and $\mathcal{T}$  if necessary). 

\begin{flushleft}
\t
({\bf Induction hypothesis}) 
Suppose that 
there is a stratum $V \simeq \R^d$ of  $\mathcal{S}$ 
such that $f$ is of class $C^1$ over the complement $\R^m - \overline{V}$. 
\end{flushleft}

Since $V$ is a $C^2$ submanifold, the tangent bundle $TV$ is of class $C^1$. 
By using the orthogonal projection to $T_xV$, 
 we can obtain a semialgebraic $C^1$ tube $T=(|T|,\pi,\rho)$ at $V$ in $\R^m$ 
as in \cite[Lemma II.5.1]{S}. 
Here $|T|$ is a semialgebraic neighborhood of $V$ in $\R^m$ 
 together with semialgebraic $C^1$ maps  $\pi: |T| \to V$ and $\rho: |T| \to \R$
so that 
\begin{itemize}
\item 
$\rho(y)\ge 0\; (y \in |T|)$ and $V=\rho^{-1}(0)$, 
\item  
$\pi \times \rho: |T|-V\to V \times\R_{>0}$ 
is a proper submersion onto $V\times(0,\,1]$, 
\item 
$d(y,V)=\inf \{||x-y||, \; x \in V\} \to 0$ ($y \in T$) 
as  $\pi(y)$ tends to a point in $\overline{V}-V$. 
\end{itemize}
We choose $T$ small enough so that 
\begin{itemize}
\item 
$|T|\cap\sigma=\emptyset$ for $\sigma\in K$ with $\sigma\cap V=\emptyset$. 
\end{itemize}

Let $D$ denote the unit closed disk  in $\R^{m-d}$ centered at $0$. 
By the construction of  $T$, 
there is a semialgebraic $C^1$ diffeomorphism 
$$g: V \times D \to |T|, \qquad g(x,u)=x+\sum_k u_k \bu_k(x)$$ 
with $x=(x_1, \cdots, x_d) \in V$ and $u=(u_1, \cdots, u_{m-d}) \in D$ 
such that 
$\bu_1(x), \cdots, \bu_{m-d}(x) \in  T_xV^\perp$ form 
a semialgebraic $C^1$ vector frame so that 
$$\pi\circ g(x,u)=x \;\;\; \mbox{and}\;\;\; 
\rho\circ g(x,u)=|u|=\left(\sum u_k^2\right)^{1/2}.$$ 
For $u\not=0$, set 
$\bar{u}:=u/|u| \in \rd D$. 

We construct a modification $\chi$ of $|T|$ as follow. 
Let $I=(0,1) \subset \R$ an open interval. 
Let $\eta:V \times I \to I$ be a semialgebraic $C^1$ map  
such that  
$I \ni t\mapsto \eta(x, t) \in I$ is a diffeomorphism for each $x \in V$, 
$\eta(x, t)=t$ for $t\ge \frac{1}{2}$, 
and 
$$\eta(x,t) \to 0 \;\; \mbox{and} \;\; 
\frac{\rd\eta}{\rd t}(x,t) \to 0 \;\;
\mbox{as}  \;\;\; t\to 0$$
(in fact, in Lemma \ref{lemma1} below, we will take $\eta$ depending on the variables $x$ 
in order to satisfy a nice property).  
We then define a semialgebraic homeomorphism $\chi: |T| \to |T|$ by 
$$
\chi(y)=\chi(g(x,u))
= g(x, \;\eta(x,|u|)\,\bar{u}) 
$$
for $y=g(x, u) \in |T|-V$ and $\chi|_V=id_V$. 
We also write $\chi(x, u)=(x, \;\eta(x, |u|)\,\bar{u})  \in V \times D$.  
Since $\chi(y)=y$ if $\rho(y) \ge \frac{1}{2}$, 
$\chi$ is extended to  
$\chi: \R^m \to \R^m$ so that $\chi=id$ outside of $|T|$.

\begin{lem}\label{lemma0}
The semialgebraic homeomorphism 
$\chi$ preserves every $\sigma \in K$ with $\sigma \cap V \not=\emptyset$, 
and $f\circ \chi: \R^m \to \R^n$ is a semialgebraic $C^0$ map and 
 is of class $C^1$ on $\R^m-\overline{V}$.  
\end{lem}

\proof 
If $V \subset \tau$ and $\sigma \cap V \not=\emptyset$, 
then $\tau$ is a face of $\sigma$ or $\tau=\sigma$; 
if $y \in |T| \cap \sigma$ and $x=\pi(y)$, 
then the segment $\overline{xy}$ is in $\sigma$, hence $\chi(y) \in \sigma$. 
Since $\pi, \rho, \eta$ are of class $C^1$, 
$\chi$ is so over $ |T|-V$.  
By the assumption that $f$ is $C^1$ off $\overline{V}$,  $f\circ \chi$ is also.  
\qed

\

\begin{lem}\label{lemma1}
We can choose $\eta(x, t)$ such that for every $1\le k \le m-d$ and for every $x \in V$, 
it holds that 
\begin{itemize}
\item 
$\displaystyle \frac{\rd (f\circ \chi)}{\rd u_k}(x,0) = 0$ 
\item 
$\displaystyle \frac{\rd (f\circ \chi)}{\rd u_k}(x',u') \to 0$,  
as  $(x', u') \to (x,0) \in V \times D$ 
\end{itemize}
\end{lem}

\proof  
For sufficiently small $t\ge 0$, we set 
$$g(t):=\sup_{x\in V, \bar{u} \in \rd D} \parallel f(x, \, t\,\bar{u})-f(x,0) \parallel.$$
Note that $f(x, \, t\,\bar{u})-f(x,0) \to 0$ 
as $x$ converges to a point in the boundary $\overline{V}-V$, 
because $f$ is continuous and $d(g(x, \, t\,\bar{u}),V) \to 0$ by a condition of $T$, 
therefore $g(t)$ is attained as the maximum. 
It is easily checked that $g(t)$ is a semialgebraic $C^0$-function; 
in particular, $g(t) \to g(0)=0$ as $t\to 0$.  
By the (ordinary) curve selection lemma, 
we see $g(t)=O(t^\alpha)$ with some positive rational number $\alpha$. 
Put $\xi(t)=t^{r_0}$ with $r_0 > \max\{\frac{1}{\alpha}, 1\}$, 
then $\frac{g (\xi(t))}{t} \to 0$ as $t\to 0$. 
Thus $\frac{1}{t}(f(x, \, \xi(t)\,\bar{u})-f(x,0))\to 0$ for any $x$ and $\bar{u}$. 
In particular,  we can choose $\eta(x,t)$ 
so that $\eta(x,t) \le \xi(t)$  for every $x \in V$; 
then every directional derivative of $f\circ \chi \; (=f(x, \eta(x,|u|)\bar{u}))$ at $(x,0)$ exists,  
and it is actually zero. Thus the first claim is verified. 
Next, we find such a function $\eta(x,t)$ which also leads to the second claim. 
Since $f$ and $\chi$ are $C^1$ on $T-V$, $f\circ \chi$ is $C^1$ and 
there are semialgebraic $C^0$ maps 
$A_k(x,u)$ and $B_k(x,u)$ for each $k$ such that 
$$
\frac{\rd (f\circ \chi)}{\rd u_k}(x,u)= \eta(x,|u|)A_k(x,u)+\frac{\rd\eta}{\rd  t}(x,|u|)B_k(x,u) 
\qquad (x \in V,\; |u|\not=0).
$$
For $x \in V$ and $0<t< 1$, set 
$$\textstyle 
a(x, t):=\min   
\left\{\; 
1, \; \frac{1}{|A_k(x,u)|}, \,\frac{1}{|B_k(x,u)|} \; 
\big| \; |u|=t, \; 1 \le k \le m-d \; \right\}.$$ 
Notice that $a(x,t)$ is a positive semialgebraic $C^0$ function on $V \times I$.   
By using a partition of unity on $V$ 
(which is locally semialgebraic in our sense), 
we can construct $\eta(x,t)$ which satisfies  
$\eta(x,t) \le \xi(t)$,  
$0<\eta(x,t) <  t \cdot a(x, t)$ and $0<\frac{\rd \eta}{\rd t}(x,t) <  t \cdot a(x, t)$ 
for $0<t\ll 1$; 
then both $\eta\cdot |A_k|$ and $\frac{\rd\eta}{\rd  t} \cdot |B_k|$ are less than $t=|u|$, 
thus we have $|\frac{\rd (f\circ \chi)}{\rd u_k}(x,u)|< 2|u|$, 
that implies the second claim.  
In fact, it is easy to find such a function $\eta(x,t)$ locally. 
For a compact semialgebraic subset $K \subset V$, 
which is the closure of an open set of $V$, 
put $a_K(t):=\min_{x \in K} a(x,t)$;  
by the curve selection lemma again, 
we may write $a_K(t)=O(t^\alpha)$ with some $\alpha\ge 0$; 
then we can take $\eta_K(t)= t^r$ ($r>r_0$) so that 
$\eta_K$ and $\frac{d}{dt}\eta_K$ are smaller than $t \cdot a_K(t)$ near $t=0$. 
By a semialgebraic bump function, 
it is extended to a semialgebraic diffeomorphism $\eta_K: I \to I$ 
with $\eta_K(t)=t$ for $t\ge 1/2$. 
Put $\eta_K(x,t):=\eta_K(t): K \times I \to I$. 
To obtain the desired function $\eta(x,t)$ on $V \times I$, 
it suffices to gather together such $\eta_K$  
by mean of a semialgebraic partition of unity on $V$. 
\qed

\

Here is an obvious claim in elementary calculus, which will be used in the proof of the next Lemma \ref{lemma2}. 

\begin{lem}\label{lemma-1}
Let $\nu(s,t)$ be a $C^1$ function on $\R\times (0,1)$ such that 
$|\frac{\rd \nu}{\rd s}(s,t)| > K$ for some fixed $K>0$. 
If $\nu$ extends continuously to $\R \times \{0\}$, 
then $\nu(s, 0)$ is not constant. 
\end{lem}

\proof Assume $\frac{\rd \nu}{\rd s}(s,t) > K$. 
Use $\nu(s_1, t) - \nu(s_0, t)> K (s_1-s_0)$ if $s_1 > s_0$, 
and take $t \to 0$. 
\qed 

\

\begin{lem}\label{lemma2} 
Let $\eta$ be as in Lemma \ref{lemma1}. 
Then there is a semialgebraic subset $S$ of $V$ with 
$\dim S < d=\dim V$ satisfying that 
\begin{itemize}
\item $\displaystyle \frac{\rd (f\circ \chi)}{\rd x_j}(x,0)=\frac{\rd f}{\rd x_j}(x,0)$ 
\item 
$\displaystyle \frac{\rd (f\circ \chi)}{\rd x_j}(x',u') \to \frac{\rd f}{\rd x_j}(x,0)$, 
as  $(x', u') \to (x,0) \in (V-S) \times D$ 
\end{itemize}
for every $1\le j \le d$ and for every $x \in V-S$. 
\end{lem}

\proof 
The first claim follows from the condition that $\chi|_V=id_V$. 
For the second, we set 
$$F(x,u)= f\circ \chi(x,u)-f(x,0), \quad F=(F_1, \cdots , F_n)  \in \R^n,$$
which is continuous on $V\times D$ and of class $C^1$ on $V\times (D-\{0\})$. 
Note that 
$$F(x,0)=0 \qquad (x \in V).$$ 
Put 
$$h(x,u):= \max_{i, j} \left|\; \frac{\rd F_i}{\rd x_j}(x,u)\; \right| \qquad (u\not=0),$$ 
which is a semialgebraic $C^0$ function on $V\times (D-\{0\})$. 
It holds that 
$$
\left|\;\frac{\rd (f_i\circ \chi)}{\rd x_j}(x',u') - \frac{\rd f_i}{\rd x_j}(x,0)\;\right| 
\le h(x',u')+\left|\;\frac{\rd f_i}{\rd x_j}(x',0) - \frac{\rd f_i}{\rd x_j}(x,0)\;\right|$$ 
and $\frac{\rd f_i}{\rd x_j}(x,0)$ is continuous by the assumption 
that $f|_V$ is $C^1$, hence 
it is enough to consider the condition that $h(x',u') \to 0$ as $(x',u') \to (x,0)$. 

Let $S$ be the subset  of $V$ 
consisting of $x$ such that 
there exists a sequence of points $(x',u') \in V \times D$ converging to $(x,0)$ 
with $h(x',u')$ not converging to $0$, that is, 
$$S=\left\{x \in V \; \bigg| 
\; \exists \, \epsilon>0\;\; \forall \, \delta>0 \;\; \exists (x', u') \in V \times D
 \;\; \mbox{s.t.}\; 
 \begin{array}{c}
 |x'-x|+|u'| < \delta,\\ h(x',u') > \epsilon
 \end{array}\; \right\}.$$
The `bad set' $S$ is semialgebraic; 
indeed, $S$ is obtained by projecting a semialgebraic subset $S'$ in $V \times \R_{>0}$ 
consisting of $(x, \epsilon)$ which satisfy the corresponding condition 
({\it i.e.},  $\forall\, \delta>0, \; \exists\, (x', u')$ such that 
$|x'-x|+|u'| < \delta$ and $h(x',u') > \epsilon$). 
We show that $S$ is of dimension $<d=\dim V$.

Suppose that $\dim S=d$. Below we make five steps for deducing the contradiction. 
A major tool is the curve selection lemma (Lemma \ref{csl}), 
which is used in each step from (1) to (4).

\

\t
(1) 
There is a non-empty semialgebraic open subset $U_0 \subset S$ 
with a constant $\epsilon_1>0$ 
such that  for each $x \in U_0$, there is a sequence  
of points $(x',u') \in V \times D$ converging to $(x,0)$ 
with $h(x', u') > \epsilon_1$. 
This is easily verified by applying Lemma \ref{csl} to 
the above semialgebraic set $S' \subset S \times \R$; 
in fact, the closure $\overline{S'}$ contains $S \times 0$, hence 
there are $U_0$ and $\epsilon_1>0$ so that 
$U_0 \times (0, 2\epsilon_1) \subset S'$; 
in particular,  $U_0 \times \epsilon_1 \subset  S'$. 
 
 \
 
 \t
 (2) 
 Let $Z:=\{(x, u) \in V \times D \; | \; h(x, u) > \epsilon_1\}$, 
then $\overline{Z}$ contains $U_0 \times \{0\}$.  
By  Lemma \ref{csl}, 
there is a non-empty semialgebraic open subset $U_1$ of $U_0$ 
which admits  a semialgebraic $C^0$ map 
$$c:U_1\times [0, a_1) \to  D$$
for some $a_1>0$ so that 
$u=c(x,t)$ is $C^1$ for $t>0$, and 
$c(x, 0)=0$, $t=|c(x,t)|$ and 
$$ h(x, c(x,t))>\epsilon_1 
\quad (\, \forall (x, t) \in  U_1 \times (0, a_1)\,).$$

\

\t 
(3) Put $c(x,t)=(c_1(x,t), \cdots, c_{m-d}(x,t)) \in D$. 
We claim that all the derivatives $\frac{\rd c_k}{\rd x_j}(x, t)$
are bounded 
over some open subset $U \times (0, a)$ of $U_1 \times (0, a_1)$. 
Suppose that $v(x,t):=|\frac{\rd c_1}{\rd x_1}(x, t)|$  is unbounded 
over any open subsets of the form $U \times (0, a)$. 
Take $M>0$. Then  for any $x \in U_1$, 
there is a sequence of points $(x', u')$ converging to $(x,0)$ 
with $v(x', u') > M$. 
Let 
$$Z':=\{(x, u) \in U_1 \times (0, a_1) \; | \; v(x, u) > M\},$$ 
then  $\overline{Z'}$ contains $U_1 \times \{0\}$, and 
$Z'$ is semialgebraic, for $v(x,t)$ is semialgebraic. 
By  Lemma \ref{csl}, 
there is an open subset $U \times (0, a)$,  over which $v(x,u) > M$. 
Take a small segment in $U$ given by  
$\gamma(s):=b+(s,0, \cdots, 0)\;\; (0\le s \le \delta, \; b \in U)$ 
and put $\nu(s,t):=c_1(\gamma(s),t)$. 
Since $c(x,0)=0$, we have $\nu(s,0)=0$, but it makes the contradiction 
by Lemma \ref{lemma-1}. 
Repeating the same argument for all $j$ and $k$, the claim is proved. 

\

\t
(4) 
Without loss of generality,  we may assume that 
the partial derivative of the first component function of $F$ 
with respect to the first coordinate of $V$ satisfies 
$$\left|\frac{\rd F_1}{\rd x_1}(x, c(x,t))\right|>\epsilon_1 
\quad (\, \forall (x, t) \in  U \times (0, a)\, ).$$
In fact, if not so, then for each $i, j$, it holds that  for any $x \in U$, 
there is a sequence $(x', t') \to (x, 0)$ 
with $\left|\frac{\rd F_i}{\rd x_j}(x', c(x', t'))\right|\le \epsilon_1$. 
Hence by applying Lemma \ref{csl} again, 
we find a (common) small positive number $a'$ 
so that 
$\left|\frac{\rd F_i}{\rd x_j}(x, c(x,t))\right|\le \epsilon_1$ 
on $U\times (0, a')$ for every $i, j$, 
but that contradicts the condition $h> \epsilon_1$ in the step (2). 
 
\
 
\t
(5) 
Consider a $C^0$ function on 
 $U\times [0, a)$ given by 
$$\phi(x,t):=F_1(x, c(x,t)),$$
which is $C^1$ off $t=0$. By the chain rule, we have  
$$
\frac{\rd \phi}{\rd x_1}(x,t)
=
\frac{\rd F_1}{\rd x_1}(x, c(x,t))
+\sum_k \;\frac{\rd (f_1\circ \chi)}{\rd u_k}(x,c(x,t))\;\frac{\rd c_k}{\rd x_1}(x,t). 
$$
Since $\left|\frac{\rd c_k}{\rd x_j}(x, t)\right|$ is bounded by (3) above, 
the sum $\sum_k$ in the right hand side converges to $0$ as $t\to 0$ by Lemma \ref{lemma1}. 
Therefore, by (4), 
we can take an open subset $U'$ with $\overline{U'} \subset U$, 
and $0<a'< a$ such that 
$$\left|\frac{\rd \phi}{\rd x_1}(x,t)\right|>\frac{\epsilon_1}{2}
\quad (\, \forall (x, t) \in  U' \times (0, a')\, ).$$
Take a segment in $U'$ of the form $\gamma(s):=b+(s,0, \cdots, 0)$ 
and put $\nu(s,t):=\phi(\gamma(s),t)$. 
By Lemma \ref{lemma-1}, $\nu(s,0)$ is not constant. 
However, for any $x \in U$, 
$$\phi(x,0)=F_1(x, c(x,0))=F_1(x, 0)=0, $$
that causes the contradiction. 
This completes the proof.\qed

\

\t 
{\it Proof of Theorem \ref{thm2}}: 
From the above Lemmas \ref{lemma1} and \ref{lemma2}, 
 $$d(f\circ \chi)(y) \to d(f\circ \pi)(x)\;\;\; \mbox{as} \;\;\; 
|T| \ni y \to x\in V-S,$$
thus $f\circ \chi$ is of class $C^1$ 
on the complement to  $(\overline{V}-V)\sqcup S$.  
The bad subset is semialgebraic and of dimension less than $d=\dim V$. 
By Lemma \ref{lemma0}, $\chi$ preserves $K$. 
By applying  Lemma \ref{C1_stratif} again we obtain 
a stratification of $f\circ \chi: \R^m \to \R^n$ 
which subdivides $K$ and the bad subset, 
hence the induction step is verified. 
This completes the proof of Theorem \ref{thm2}. \qed

\begin{rem}\label{PA}{\rm 
({\bf PA spaces}) 
In Theorems \ref{thm1} 
$X$ is assumed to be embedded in $\R^m$, 
but the global embeddability  is not needed essentially, because   
key lemmas in the proof deal with local problems. 
Indeed, our argument fits with the framework of 
(possibly non-compact) {\it piecewise algebraic spaces} 
 introduced by Kontsevich-Soibelman \cite{KS}; 
a PA space is a locally compact Hausdorff space $X=\bigcup X_\alpha$ 
made by glueing countably many compact semialgebraic sets $X_\alpha$ 
via semialgebraic homeomorphisms on constructible subsets 
$X_\alpha \supset X_{\alpha\beta} \stackrel{\simeq}{\to} X_{\beta\alpha}\subset X_\beta$. 
In the same approach, we can see that 
any PA space admits a (PA) triangulation 
 with a $C^1$ realization so that 
the (locally finite) simplicial complex is locally realized in some $\R^m$. 
}
\end{rem}

\section{Semialgebraic chains and integrations}

\subsection{Integrals of differential forms}
Let $M$ be a semialgebraic manifold,
 $X$ a compact semialgebraic subset of $M$ ($\dim X=p$) 
 with a semialgebraic subset $\rd X \subset X$ so that  
 $\rd X$ is locally  the boundary of the regular part $X_{\rm reg}$ 
 around each regular point of $\rd X$. 
Assume that $X$ has a fundamental class $\mu_{X} \in H_p(X, \rd X; \Z)$, 
{\it i.e.},  $i_*\mu_{X}$ generates $H_p(X, X-x; \Z)\simeq \Z$ 
for any regular point $x \in X _{\rm reg}-\rd X$ with inclusion $i: X-x \hookrightarrow X$. 

Let $(K, f)$ be a semialgebraic triangulation of $X$ with a $C^1$ realization 
as in Theorem \ref{thm1}. 
Then for each simplex $\sigma \in K$, 
$f|_\sigma: \sigma \to X$ is a semialgebraic $C^1$ map. 
The  chain $\sum \sigma$ of all $p$-simplices in $K$ with suitable orientation 
represents $\mu_X$; we  define the integral over $X$ 
of a differential $p$-form $\omega$ on $M$ by 
$$\int_X \omega := \sum \int_\sigma (f|_\sigma)^* \omega.$$ 
We show that the integral does not depend on the choice of $(K,f)$. 

Let $(K_1, f_1)$ and $(K_2,f_2)$ be 
semialgebraic triangulations of $X$ with $C^1$ realizations. 
Obviously, there is a semialgebraic homeomorphism $\tau: |K_1| \to |K_2|$ 
with $f_1=f_2\circ \tau$. 
Applying Lemma \ref{C1_stratif}, we have a stratification for $\tau$, 
and let $\mathcal{A}$ be the set of 
strata of the source stratification of $\tau$, $\sigma \in K_1$ 
and $\tau^{-1}(\sigma')$ of $\sigma' \in K_2$. 
By Theorem \ref{triangulation}, 
we can find a semialgebraic triangulation $(K, \phi_1)$ of $|K_1|$ compatible with $\mathcal{A}$. 
Set $\phi_2=\tau \circ\phi_1:|K| \to |K_2|$ and 
 $g=f_1 \circ \phi_1=f_2 \circ \phi_2: |K| \to X$. 
Notice that any simplex of $K_j \; (j=1,2)$ is 
the union of the image of finitely many simplices of $K$ via $\phi_j$. 
By Corollary \ref{thm2_cor}, 
there is a common panel beating $\chi$ of $|K|$ 
such that $\phi_j \circ \chi\; (j=1,2)$ are of class $C^1$. 
Put $f:=g\circ \chi=f_j\circ (\phi_j\circ \chi)$, then $f$ is $C^1$,  
for $f_j$ and $\phi_j\circ \chi$ are so. 
We call the obtained triangulation $(K, f)$ of $X$ 
a  common $C^1$ refinement of $(K_1, f_1)$ and $(K_2,f_2)$. 

$$
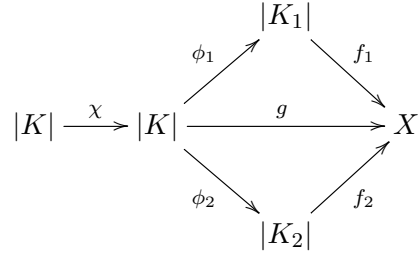
\xymatrix{
&& |K_1| \ar[dr]^{f_1} & \\
|K| \ar[r]^\chi & |K| \ar[ur]^{\phi_1} \ar[dr]_{\phi_2} \ar[rr]^g && X\\
&& |K_2| \ar[ur]_{f_2} & \\
}
$$

\begin{thm}\label{thm3}
The integral $\int_X \omega$ is well-defined, {\it i.e.},  
it does not depend on the choice of 
semialgebraic triangulation of $X$ with $C^1$ realization. 
\end{thm}

\proof 
Let $\sigma: \Delta^p \to X$ be a semialgebraic $C^0$ map. 
By a panel beating, we may assume that $\sigma$ itself is of class $C^1$. 
Let $(K, h)$ be a semialgebraic triangulation of $\Delta^p$ with a $C^1$ realization. 
Then $h: \Delta^p=|K| \to \Delta^p$ is a $C^1$ map 
which is diffeomorphic off a measure zero subset, 
thus we have 
$\int_\sigma \omega=\int_{\Delta} g(x) dx = \int_{\Delta} g(h(u)) |\frac{\rd h}{\rd u}| du=\int_{\sigma\circ h}\omega$ for a differential $p$-form $\sigma^*\omega=g(x)dx_1\wedge \cdots \wedge dx_p$. 
By using a common $C^1$ refinements of triangulations, the assertion immediately follows. 
 \qed

\begin{prop}
The Stokes formula holds, $\int_X d\omega= \int_{\rd X} \omega$, provided 
$X$ and $\rd X$ have fundamental cycles. 
\end{prop}

\proof 
It is the same as in the smooth case. 
\qed

\

The above argument is also applied to 
the case of non-compact $X$ and compact supported $\omega$. 
We can also work on  subanalytic sets or more generally $\frak X$-sets 
(see Remarks below); 
in particular, 
the same treatment is available for complex analytic sets and 
integrations of $(p,q)$-differential forms. 
There is an alternative approach for the integration 
on subanalytic sets, e.g. \cite{Valette15}.

\begin{rem}\label{current}{\rm 
({\bf Semialgebraic currents}) 
Our integral $\int_X \omega$ over a compact semialgebraic set $X$ 
coincides with the integral as a {\it semialgebraic current}  
in geometric measure theory \cite{Federer, Hardt, HLTV}. 
The key point in the theory is that $X$ has a finite volume, 
hence the integral over the interior absolutely converges. 
Subanalytic currents are also available. 
However,  
the theory of currents is unclear 
for a more general setting 
such as the category of $\frak X$-sets 
or o-minimal category over any real closed field, 
while our approach does not have any trouble in such a general context. 
}
\end{rem}

\begin{rem}\label{CDGA}{\rm 
({\bf de Rham homotopy theory of PA spaces}) 
In \cite{KS} 
the foundation of {\it semialgebraic} differential forms on 
semialgebraic sets, more precisely {\it PA forms} on PA spaces, has been introduced 
for the use in proving the formality of little cubes operads; 
in particular, as an analogy to the de Rham algebra of differential forms on smooth manifolds, 
the algebra of PA forms 
encodes the information of the real homotopy type of the PA space. 
A large portion of the theory of PA forms has been clarified in \cite{HLTV} 
through a bit long discussion using semialgebraic currents. 
Our approach using panel beatings simplifies the arguments 
in a straightforward way; it enables us to directly connect  
the de Rham algebra of PA forms with Sullivan's model (CDGA)  with real coefficients. 
Another advantage is that our approach works over an {\it arbitrary} real closed field. 
}
\end{rem}

\begin{rem}\label{o_minimal}
{\rm ({\bf O-minimal category}) 
We may work over 
any o-minimal category over any general real closed field $R$. 
Let $X$ be a bounded definable subset $X$ of dimension $d$ in $R^m$. 
Assume that 
$X \subset I^m$ 
where $I$ is the closed interval $[0,L]$, $L \in R$. 
By subdividing $I$ into $n$ intervals $[L/n, (k+1)L/n]$, 
$I^m$ breaks into $m$-dimensional small cubes $U_i$ 
 with $1\le i \le n^m$. 
For each $i$, set $c_i=1$ if $X \cap U_i \not=\emptyset$, 
$c_i=0$ otherwise. Let $\delta_n:=\sqrt{m}\cdot L/n$, the diameter of the small cube, 
and put $v_n^{(d)}(X):= \left(\sum_{i=1}^{n^m} c_i\right) (\delta_n)^d$. 
Then, by using Theorem \ref{thm1}, we can show that 
 $v_n$ is uniformly bounded in $R$, that is, 
there is $c \in R$ such that $v_n^{(d)}(X) <c$ for any $n \in \N$. 
Moreover, the {\it $d$-th  Hausdorff measure} ${\mathcal{H}}^d(X)$ 
(or the {\it volume}) of $X$ is given by $\lim_{n \to \infty} v_n^{(d)}(X)$, 
which exists in some saturated field containing $R$. 
However, the argument about the saturated field is quite complicated, 
so we do not discuss it here. See \cite{MS} for details. 
}
\end{rem}


\small\begin{thebibliography}{99}
\bibitem{BCR} J.~Bochnak, M.~Coste and M.F~Roy, 
\textit{Real algebraic geometry}, Springer-Verlag (1998).  
\bibitem{CP} M. Czapla and W. Paw\l ucki, 
Strict $C^1$-triangulations in o-minimal structures, preprint (2017). 
\bibitem{Federer} H.~Federer, \textit{Geometric measure theory}, Springer-Verlag, 
Heidelberg and New York, 1969.
\bibitem{G} C.~G.~Gibson, K.~Wirthm\"uller, A. A.~Du Plessis, E. J. N.~Looijenga, Topological stability of smooth mappings, Lecture Notes in Math., 552, Springer, 1976.
\bibitem{Hardt} R.M.~Hardt, Topological properties of subanalytic sets, Trans. Amer. Math. Soc., 211 (1975),  57--70.
\bibitem{HLTV} R.M.~Hardt, P.~Lambrechts, V.~Turchin and I.~Voli\'c, 
Real homotopy theory of semi-algebraic sets, 
Algebraic \& Geometric Topology {\bf 11} (2011), 2477--2545, DOI: 10.2140/agt.2011.11.2477. 
\bibitem{HKT} M. Hanamura, K. Kimura and T. Terasoma, 
Integrals of logarithmic forms on semi-algebraic sets and a generalized Cauchy formula, Part I: 
convergefnce theorems, 
preprint (2015), arXiv:1509.06950. 
\bibitem{KS} M. Kontsevich and Y. Soibelman, 
Deformations of algebras over operads and the Deligne conjecture, 
In Conf\'erence Mosh\'e Flato 1999, Vol. I (Dijon), 
Math. Phys. Stud. {\bf 21}, Kluwer Acad. Publ., (2000), 255--307.
\bibitem{L} S.~{\L}ojasiewicz, Triangulation of semi-analytic sets, Ann. Scu. Norm. Sup. Pisa {\bf 18} (1964), 449--474.  
\bibitem{Loi} T. L. ~Loi, 
Verdier and strict Thom stratifications in o-minimal structures, Illinois J. Math. {\bf 42} (1998), 
347--356.  
\bibitem{MS} J.~Ma\v{r}\'ikov\'a and M.~Shiota, 
Measuring definable sets in o-minimal fields,  Israel Journal of Math. {\bf 209} (2015), 687--714. 
\bibitem{MT} C.~Murolo and D.~Trotman, 
Horizontally-$C^1$ controlled stratified maps and Thom's first isotopy theorem, 
Comptes Rendus Acad. Sci. Ser. I - Math. {\bf 330} (8) (2000), 707--712. 
\bibitem{S} M.~Shiota, \textit{Geometry of subanalytic and semialgebraic sets}, 
Progress in Math. {\bf 150}, Birk\-h\"auser, 1997.
\bibitem{Suwa} 
 T.~Suwa, Indices of vector fields and residues of singular holomorphic foliations , Actualit\'es Math\'ematiques, Hermann, Paris (1998).
\bibitem{Valette15} G.~Valette, 
Stokes' formula for stratified forms, Ann. Polon. Math. {\bf 114} (2015), 197--206.
\end{thebibliography}
\end{document}